\documentclass[12pt]{amsart}

\usepackage{fullpage}
\usepackage{graphicx}
\usepackage{amsfonts}
\usepackage{amstext}
\usepackage{amssymb}
\usepackage{amsmath}
\usepackage{amsthm}
\usepackage{amsopn}
\usepackage{bm}
\usepackage[ps,dvips,matrix,curve,frame,arrow,rotate,line]{xy}

\newtheorem{theorem}{Theorem}
\newtheorem{prop}[theorem]{Proposition}
\newtheorem{lemma}[theorem]{Lemma}
\newtheorem{corollary}[theorem]{Corollary}

\theoremstyle{remark}
\theoremstyle{definition}
\newtheorem{remark}[theorem]{Remark}
\newtheorem{example}[theorem]{Example}
\newtheorem{defn}[theorem]{Definition}

\newcommand{\divides}{\mid}

\newcommand{\C}{\mathcal{C}}
\newcommand{\B}{\mathcal{B}}
\newcommand{\F}{\mathcal{F}}
\newcommand{\D}{\mathcal{D}}

\newcommand{\QQ}{\mathbb{Q}}
\newcommand{\RR}{\mathbb{R}}
\newcommand{\ZZ}{\mathbb{Z}}
\newcommand{\NN}{\mathbb{N}}
\newcommand{\CC}{\mathbb{C}}
\newcommand{\FF}{\mathbb{F}}
\newcommand{\isom}{\cong}

\newcommand{\GL}{\operatorname{GL}}

\newcommand{\ii}{\sqrt{-1}}

\newcommand{\oh}[1]{\widehat{#1}}
\newcommand{\sq}[1]{\widetilde{#1}}

\DeclareMathOperator{\Frob}{Frob}

\DeclareMathOperator{\Mat}{Mat}
\DeclareMathOperator{\Hom}{Hom}

\DeclareMathOperator{\gpid}{e}
\DeclareMathOperator{\idgp}{\gpid}
\DeclareMathOperator{\id}{id}
\DeclareMathOperator{\im}{im}
\DeclareMathOperator{\Gal}{Gal}

\DeclareMathOperator{\Mor}{Hom}
\DeclareMathOperator{\Obj}{Ob}

\newcommand{\eqdef}{\overset{\text{def}}{=}}
\newcommand{\ilim}{\mathop{\varprojlim}\limits}

\swapnumbers
\numberwithin{theorem}{section}

\newcommand{\psm}{quotient-preserving map}
\newcommand{\MD}{\mathbf{MD}}
\newcommand{\TopGp}{\mathbf{TopGp}}

\begin{document}
\title{Measurable Dynamics of Maps on Profinite Groups}
\author{James Kingsbery}
\address[James Kingsbery]{Department of Mathematics\\
     Williams College \\ Williamstown, MA 01267, USA}
\email{James.C.Kingsbery@williams.edu}
\author{Alex Levin}
\address[Alex Levin]{Harvard University, MA 02138, USA}
\email{levin@fas.harvard.edu}
\author{Anatoly Preygel}
\address[Anatoly Preygel]{Harvard University, MA 02138, USA}
\email{preygel@fas.harvard.edu}
\author{Cesar E. Silva}
\address[Cesar E. Silva]{Department of Mathematics\\
     Williams College \\ Williamstown, MA 01267, USA}
\email{csilva@williams.edu}
\begin{abstract}
We study the measurable dynamics of transformations on profinite groups, in particular of those which factor through sufficiently many of the projection maps; these maps generalize the $1$-Lipschitz maps on $\ZZ_p$.
\end{abstract}
\date{\today}
\maketitle

\section{Introduction}\label{sec:intro}
Several authors have studied the measurable dynamics of polynomial maps 
that define Haar measure-preserving transformations on balls or spheres in the (locally compact) 
field of  $p$-adic numbers, see for example \cite{GKL}, \cite{CoelhoParry}, \cite{Anashin},
\cite{BrykSilva}.  Anashin \cite{Anashin} has studied a class of maps on 
$\ZZ_p^k$ that are $1$-Lipschitz and that he calls {\it compatible};    Anashin  stated  
that if a compatible (i.e., $1$-Lipschitz ) map is measure-preserving, then it is bijective,
and moreover it is an isometry of $\ZZ_p^k$ (under the $p$-adic metric).  It is also 
true that if it is bijective then it is measure-preserving, hence an isometry (see
\cite[Lemma~4.5 ]{BrykSilva}).  It was also shown in 
\cite{BrykSilva} that an isometry on a compact-open subset of $\QQ_p$ is never totally 
ergodic, in contrast to the real case where, for example, irrational rotations 
on the circle are totally ergodic.  
  In this paper we introduce a class of maps 
called {\it {\psm}s} that generalize the asymptotically compatible (and compatible) maps 
of Anashin and classify their measurable dynamics.    However, rather then studying 
these maps on $\ZZ_p$ we find that their natural setting is in the context of profinite groups.  
We now outline the contents of the various sections.

Section~\ref{sec:invLim} reviews inverse limits and states the basic properties of 
profinite groups that we will use.  Section~\ref{sec:mdCat} is a review of applications
of these notions, in particular of inverse limits, to the context of measurable dynamics.
In Section \ref{sec:factProj} we introduce the notion of  {\psm}s and prove the 
following theorem on the dynamics of these maps.

\begin{theorem}\label{thm:profDynamics}
Let $G$ be a second-countable profinite group, $\mu$ normalized Haar measure on $G$, and $T: G \to G$ a {\psm}. Define the \emph{finite factor set of $T$} as
\[ \F(T) = \{ N \lhd_O G : T\text{ factors through }\pi_N:G \to G/N \}. \]  Let $\F \subseteq \mathcal{F}(T)$ be a base for the neighborhoods of $\idgp \in G$.  For each $N \in \mathcal{F}(T)$ let $T_N$ denote the induced map $G/N \to G/N$. Then, the following are equivalent:
\begin{enumerate}
\item $T$ is measure-preserving (equivalently nonsingular) with respect to $\mu$;
\item $T_N$ is bijective for each $N \in \F$;
\item $T$ is surjective;
\item There exists a translation invariant metric $d$ inducing the topology on $G$ such $T$ is an isometry with respect to $d$ and the subset of $\F$ consisting of sets that are balls of some radius with respect to $d$ is a base for the neighborhoods of $\idgp \in G$.
\end{enumerate}

Also, the following are equivalent:
\begin{enumerate}
\item $T$ is measure-preserving and ergodic with respect to $\mu$;
\item $T_N$ is measure-preserving and ergodic with respect to $\mu_{G/N}$ for eaah $N \in \F$;
\item $T_N$ is minimal with respect to $\mu_{G/N}$ for each $N \in \F$.
\end{enumerate}
\end{theorem}

Section~\ref{sec:hom} applies our methods to the case of continuous homomorphisms, where the additional structure allows us to give a simpler characterization of {\psm}s.  Finally, Section~\ref{sec:prod} applies our results to products of {\psm}s.  The prototypical
examples of such products are given by products of $1$-Lipschitz maps on $\ZZ_p$, for possibly
different primes $p$.  The main result of Section~\ref{sec:prod} is Theorem~\ref{thm:prodProfErg}.

\subsection{Acknowledgements}
This paper is based on research by the Ergodic Theory group of the 2005 SMALL  summer research project at Williams College and was first posted on Silva's website on November 9, 2005.   Support for the project was provided by National Science Foundation REU Grant DMS - 0353634 and the Bronfman Science Center of Williams College.

\section{Inverse limits}\label{sec:invLim}
%We say that a pair $(I,\leq)$ is a \emph{directed set} if the following conditions are satisfied:
%\begin{enumerate}
%\item $I$ is a set, and $\leq$ is a partial order on $I$;
%\item for each $\alpha,\beta \in I$ there exists a $\gamma \in I$ such that $\alpha \leq \gamma$ and $\beta \leq \gamma$.
%\end{enumerate}
%We may regard a directed set as a category with objects the elements of $I$, and with a morphism $i \to j$ for $i \leq j$.

For our purposes, we are primarily interested in inverse limits in two categories:
\begin{enumerate}
\item The category $\TopGp$: The objects of $\TopGp$ are topological groups, and the morphisms are continuous group homomorphisms.
\item The category $\MD$: The objects are measurable dynamical systems, and the morphisms are measure-preserving maps commuting (almost everywhere) with the action of the dynamical systems (identifying two morphisms if they agree almost-everywhere).
\end{enumerate}
In the following, we let $\mathfrak{C}$ be an arbitrary category; in light of the above, the reader should feel free to replace it with either of the above.

A \emph{inverse system in $\mathfrak{C}$}, denoted $\D: (I,\leq) \to \mathfrak{C}$ consists of the following data:
\begin{enumerate}
\item A directed set $(I,\leq)$ (i.e. $(I,\leq)$ is a partially ordered set, such that each finite subset has an upper bound in $I$);
\item A collection $\{\D(i) \in \Obj_{\mathfrak{C}} : i \in I \}$ of objects of $\mathfrak{C}$;
\item A collection $\{\D(i,j) \in \Mor_{\mathfrak{C}}(\D(j),\D(i)) : i \leq j\}$ of morphisms such that for all $i \leq j \leq k \in I$ we have $\D(i,j) \circ \D(j,k) = \D(i,k)$ and such that $\D(i,i)  = \id_i$ for all $i \in I$.
\end{enumerate}

A pair $(L,\{\pi_i\})$ with $L \in \Obj_{\mathfrak{C}}$ and with $\{ \pi_i \in \Mor_{\mathfrak{C}}(L,\D(i)): i \in I\}$ a collection of morphisms such that $\D(i,j) \circ \pi_j = \pi_i$ for all $i \leq j \in I$ is said to satisfy the \emph{defining property of an inverse limit} for the inverse system $\D$.

An \emph{inverse limit} for the inverse system $\D$ is a pair $(L,\{\pi_i\})$ satisfying the defining property of an inverse limit and the following universal property: For any pair $(L',\{\pi'_i\})$ satisfying the defining property of an inverse limit there must exist a unique morphism $L' \to L$ making the following diagram commute for all $i \leq j \in I$:
\[\xymatrix{
   &  L' \ar@{-->}[dd]^{!} \ar[dddl]_{\pi'_j} \ar[dddr]^{\pi'_i} & \\ \\
   &  L  \ar[dl]^{\pi_j} \ar[dr]_{\pi_i} & \\
H_j \ar[rr]_{\D(i,j)} & & H_i } \]
Such an object, which is unique if it exists, is denoted by
\[ \ilim_{i \in I} \D(i) .\]  %As is usual for objects defined via a universal construction, if an inverse limit exists exist then it is unique up to canonical isomorphism in $\mathfrak{C}$.
% That is, say $(L,\{\pi_i\})$ and $(L',\{\pi'_i\})$ are both inverse limit objects for $\mathcal{D}$.  Then, the above definition guarantees us morphisms $\phi: L \to L'$ and $\varphi: L' \to L$ such that the appropriate diagrams commute.  Now, $\pi'_i \circ \phi \circ \varphi = \pi'_i$, so the universal property gives that $\phi \circ \varphi = \idgp_{L'}$; similarly $\varphi \circ \phi = \idgp_L$.  So, $\phi$ and $\varphi$ are two-sided inverses, giving the canonical isomorphism of the two inverse limit objects.

If $\mathfrak{C}$ is $\TopGp$ then each directed system in $\mathfrak{C}$ has an inverse limit, given by the following construction:
\[ \ilim_{i \in I} \D(i) = \{ x \in \prod_{i \in I} \D(i) : \pi_i(x) = \D(i,j)(\pi_j(x))\text{ for all $i \leq j \in I$} \} , \] with the subspace topology from the product topology and with projection maps given by the projection maps from the product.

We are now ready to define a \emph{profinite group}.  We say that a topological group $G$ is \emph{profinite} if it is isomorphic, as a topological group, to an inverse limit of finite groups.  That is, if %$G$ is profinite if and only if there exists an directed set $(I,\leq)$, and an inverse system of \emph{finite} (topological via the discrete topology) groups $\D: (I,\leq) \to \mathbf{TopGp}$ such that 
\[ G \isom \ilim_{i \in I} \D(i) \] for $\D: (I,\leq) \to \TopGp$ an inverse system of \emph{finite} (topological via the discrete topology) groups.

Let us sketch and cite some standard results on profinite groups:
\begin{prop}\label{prop:propProfGp}
Let $G$ be a profinite group.  Then:
\begin{enumerate}
\item $G$ is a compact Hausdorff totally-disconnected topological group.  Moreover, these properties characterize profinite groups.
\item Every open subgroup $U \le_O G$ is also closed (this in fact holds for all topological groups).
\item Every open subgroups $U \le_O G$ has finite index.
\item The normal open subgroups form a base for the neighborhoods of $\idgp \in G$ (equivalently, their translates form a base for the topology on $G$).
\item Let $\F$ be a collection of open normal subgroups of $G$ such that $\F$ is a base for the neighborhoods of $\idgp \in G$.  Then, we may order $\F$ by inclusion, and for $N \supseteq N'$ we have a projection $G/N' \to G/N$.  This makes the system of quotients $G/N$ into an inverse system, with
\[ G \isom \ilim_{N \in \F} G/N, \] where the inverse limit and isomorphism are $\TopGp$.
\item Let $\B$ be smallest $\sigma$-algebra containing the compact subsets of $G$.  Then, there is a unique measure $\mu$ on $\B$ such that $\mu(gS)=\mu(sG)=\mu(S)$ for $g \in G$ and $S \in \B$, $\mu$ is regular, and  $\mu(G)=1$.  We call $\mu$ the (normalized) Haar measure on $G$. \label{prop:propProfGp-vi}
\end{enumerate}
\end{prop}
\begin{proof}
For (i), note that the product space in the construction given above is compact Hausdorff.   Then, $G$ corresponds to a closed subgroup of the product, and so is also a compact Hausdorff topological group.  That $G$ is totally disconnected then follows from (ii) and (iv).  For the converse, it suffices to show that (iv) holds for such a space and then use the proof of (v); for this see the reference below.

Distinct cosets of $U$ are disjoint; so the union of the cosets different from $U$ is just $G \setminus U$, and this set must be open.  This proves claim (ii).  Claim (iii) follows by compactness.

Say $G \isom \ilim_{i \in I} \D(i)$, $\D(i)$ finite groups with the discrete topology, and let $\pi_i: G \to \D(i)$ be the projection map.  Then, $\ker \pi_i$ is a normal open subgroup of $G$ for each $i \in I$.  We readily check that these form a base for the neighborhoods of $\idgp \in G$ (indeed, their cosets are just the restriction of the standard base for the product topology on the inverse limit).  This proves (iv).

Now, say $\F$ forms a base for the neighborhoods of $\idgp \in G$.  Let $\pi_N: G \to G/N$ be the quotient maps.  Then, $(G,\{\pi_N\})$ satisfies the defining property of the inverse limit, so by the universal property of the inverse limit we have a canonical map
\[ \phi: G \longrightarrow \ilim_{N \in \F} G/N \]
such that the appropriate diagram must commute.  Note that this map must be an injection, for \[ \bigcap_{N \in \F} \ker \pi_N = \bigcap_{N \in \F} N = \{ 1 \}. \]  Furthermore, the image of $\phi$ must be dense, and must be compact as $G$ is compact, $\phi$ continuous, and the inverse limit Hausdorff.  So, $\phi$ is surjective.  So, $\phi$ is a continuous bijection.  But, $\phi$ must take closed, hence compact, sets to compact, hence closed, sets; so $\phi^{-1}$ is continuous.  So, $\phi$ is an isomorphism of topological groups.  This proves (v).  For more on the general theory of topological groups see for instance \cite{Pontryagin}.  For complete proofs of the above claims, see for instance \cite[p. 17-20]{Wilson}.

Finally, $G$ is a compact topological group, so it is unimodular and has a unique (left and right) Haar measure.  This proves (vi).  For more details on Haar measure on locally compact groups and the unimodularity of compact groups see for instance \cite[p. 36-47]{Folland}.
\end{proof}

\begin{example}
Let $I = \NN$, and for $k \in I$ let $\D(i) = \ZZ/p^i \ZZ$.  For $i \leq j \in I$ let $\D(i,j): \ZZ/p^j \ZZ \to \ZZ/p^i \ZZ$ be the reduction $\bmod p^i$ map.  Then, we have
\[ \ZZ_p \isom \ilim_{i \in I} \D(i) = \ilim_{k \geq 1} \ZZ/p^k \ZZ, \] where $\ZZ_p$ refers to the additive group of the ring of $p$-adic integers.
\end{example}

\section{Measurable dynamical structure}\label{sec:mdCat}
By a \emph{measurable dynamical system} we mean a $4$-tuple $(X,\mu,\B,T)$ where $X$ is a set, $\B$ is a $\sigma$-algebra of subsets of $X$, $\mu$ is a probability measure on $\B$, and $T$ is a $\B$-measurable function.  We define a \emph{morphism of measurable dynamical systems} $(X,\mu,\B,T) \to (X',\mu',\B',T')$ to be an equivalence class of maps $\phi: X \to X'$ such that $\phi$ is measurable and measure-preserving, and $\phi \circ T(x) = T' \circ \phi(x)$ holds outside a set of $\mu$-measure $0$; our equivalence relation is to identity $\phi: X \to X'$ and $\phi': X \to X'$ when $\phi(x)=\phi'(x)$ holds outside a set of $\mu$-measure $0$.  These definitions define a category, which we shall denote $\MD$.

Inverse limits need not always exist in $\MD$; indeed even when the inverse system consists just of finite direct products, there need not be a measure on the topological inverse limit \cite[p. 214]{Halmos}.  There are significant existence results, such as in the case of standard spaces \cite{Parthasarathy} or of topological measures on compact spaces \cite{Choksi}.  Even without these topological restrictions, we may sometimes be guaranteed that an inverse system has an inverse limit; furthermore when an inverse limit exists its dynamics are closely related to the dynamics of the systems in the inverse system:
\begin{prop}\label{prop:propMDIlim}
Let $(I,\leq)$ be an directed set, and $\D: I \to \MD$ an inverse system in $\MD$.  Moreover, assume there is an object $U=(X,\mu,\B,T)$ and morphisms $\{ \pi_i \in \Mor_{\MD}(U \to \D(i)) : i \in I\}$ such that the following diagram commutes for each $i \leq j \in I$
\[ \xymatrix{ & U \ar[dl]_{\pi_j} \ar[dr]^{\pi_i}\\ \D(j) \ar[rr]_{\D(i,j)} & &  \D(i)} \]  For each $i \in I$, let $\B_i$ denote the $\sigma$-algebra of measurable sets of $\D(i)$, and let $\sq{\B}$ be the smallest $\sigma$-algebra containing
\[ \bigcup_{i \in I} \pi_{i}^{-1}(\B_i). \]  Then, $(X,\mu,\sq{\B},T)$ is an inverse limit for $\D$.

Moreover, if $L$ is an inverse limit for $\D$ then $L$ is measure-preserving if and only if $\D(i)$ is measure-preserving for each $i \in I$.  The previous sentence still holds when one adds to ``measure-preserving'' any of the following additional conditions: ergodic, weakly mixing, mixing.
\end{prop}
\begin{proof}
See \cite{Brown}.
\end{proof}

Now, Proposition~\ref{prop:propProfGp}(\ref{prop:propProfGp-vi}) turns each profinite group, in a natural way, into a probability space.  Say $G$ is a profinite group, $\mu$ Haar measure on $G$, and $\B$ the $\sigma$-algebra of $\B$-measurable sets.  Then, for any $\mu$-measurable map $T: G \to G$ we have that the $4$-tuple $\Sigma = (G,\mu,\B,T)$ is an object of $\MD$.  The final statement of Proposition~\ref{prop:propProfGp} combined with Proposition~\ref{prop:propMDIlim} suggests that we may be able to study the dynamics of a system on $G$ by looking at systems on some finite quotients of $G$.  Unfortunately, for $N \lhd_O G$ an open normal subgroup, $T$ need not induce a well-defined map $G/N \to G/N$.  We may recover some such information through the following construction.

For $N \lhd_O G$ we define the following objects:
\begin{itemize}
\item Let \[ X_N = \prod_{k \geq 0} G/N, \] let $\pi_N: G \to G/N$ be the quotient map, and let the map $\Phi_N: G \to X_N$ be given by \[ x \mapsto (\pi_N(x), \pi_N(Tx), \pi_N(T^2 x), \pi(T^3 x), \ldots) \qquad \text{that is $\varpi_k \circ \Phi = \pi_N \circ T^k$}, \] where $\varpi_k: X_N \to G/N$ is projection to the $k^\text{th}$ slot. 
\item We may define a measure on $X_N$ such that $\Phi_N$ is measure-preserving; specifically, let $\mu_N = \mu \circ \Phi_N^{-1}$, let $\B_N$ the $\sigma$-algebra of $\mu_N$-measurable sets.
\item  Finally, let $T_N$ be the left-shift map on $X_N$.  Then, we may define the following measurable dynamical system:
\[ \Sigma_N = (X_N,\mu_N,\B_N,T_N). \]
\end{itemize}

\begin{lemma}\label{lem:ilimProf}
Let $\Sigma = (G,\mu,\B,T)$ be a measurable dynamical system with $G$ a profinite group and $\mu$ Haar measure on $G$.  Let $\Sigma_N, \Phi_N$ be as above.

Say $I \subseteq \{ N \lhd_O G \}$ is ordered by set-inclusion.  For $N \supseteq N' \in I$, we have a natural projection $G/N' \to G/N$; this induces a morphism (of $\MD$) $\Sigma_{N'} \to \Sigma_N$.  Now, we may define $\D: (I,\supseteq) \to \MD$ by
\[ \D(N) = \Sigma_N \qquad \D(N,N')=\text{the above morphism $\Sigma_{N'} \to \Sigma_{N}$} \]
for all $N,N' \in I$.

Then:
\begin{enumerate}
\item $\D$ is an inverse system in $\MD$;
\item $(\Sigma,\{ \Phi_N \})$ satisfies the defining property for the inverse limit of $\D$;
\item $\D$ has an inverse limit in $\MD$;
\item If $G$ is second-countable and $I$ forms a base for the neighborhoods of $\idgp \in G$, then $(\Sigma,\{ \Phi_N \})$ is an inverse limit for $\D$.
\end{enumerate}
\end{lemma}
\begin{proof}
The commutativity of the appropriate diagrams for (i) and (ii) are routine verifications.  We note that the maps $\pi_N$, as well as the maps $\D(N,N')$ are surjective continuous group homomorphisms.  It is a standard result that surjective continuous group homomorphisms preserve Haar measure.  Also, for each $N \in I$, the map $\Phi_N$ is continuous and is measure-preserving by construction of $\mu_N$.  So, all relevant maps are indeed morphisms in $\MD$ and claims (i) and (ii) are complete.  Then, claim (iii) follows by Prop~\ref{prop:propMDIlim}.

Now, by Prop~\ref{prop:propMDIlim}, letting $\sq{\B}$ be the smallest $\sigma$-algebra containing
\[ \bigcup_{N \in I} \Phi_N^{-1}(\B_N), \] we have that $(X,\mu,\sq{\B},T)$ is an inverse limit for $\D$.  Noting that the maps $\Phi_N$ are measurable we have $\sq{\B} \subseteq \B$.

Say $G$ is second-countable.  Each element of $I$ is a compact-open set, and is thus a finite union of elements of the countable base of $G$.  As the collection of finite subsets of a countable set is itself countable, we have that $I$ must be at most countable.  Moreover, each $N \subseteq I$ has finitely many distinct translates.  So, if $I$ forms a base for the neighborhoods of $\idgp \in G$, then the collection of translates of the elements of $I$ form a countable base for the topology of $G$.

For $N \subseteq I$, the cosets of $N$ are contained in $\Phi_N^{-1}(\B_N)$.  So, $\sq{\B}$ contains all translates of $I$, and hence a countable base for the open sets of $G$.  By countable unions, $\sq{\B}$ contains the open sets of $G$, and by taking complements it contains the closed sets of $G$ and so the compact sets.  Recalling that $\B$ was generated by the compact sets, we have $\B \subseteq \sq{\B}$.  With the above, this implies that $\B = \sq{\B}$ and proves our claim.
\end{proof}

\begin{example}
Let $G=\ZZ_p$.  Note that each element of $\ZZ_p$ has a unique expression of the form $c + p d$ with $c \in \{ 0, \ldots, p-1\}$ and $d \in \ZZ_p$.  Then, we may define $T: G \to G$ by
\[ T(c + p d) = d \text{ for $c \in \{0,\ldots, p-1\}$, $d \in \ZZ_p$}. \]

Then, $T$ is a surjective, $p$-to-$1$, measure-preserving map.  Take $N=p \ZZ_p$.  Then, $\Sigma_N$ is a Bernoulli shift on $p$ symbols.  Moreover, one can show that the map $\Phi_N: G \to X_N$ is a measurable (and topological) isomorphism.
\end{example}

\begin{example}
Let $G = \ZZ_p$, and define the transformation $f: G \to G$ by \[ f(x) = {x \choose p} = \frac{x(x-1)\cdots(x-p+1)}{p!}. \]  

Take $N = p \ZZ_p$.  It is possible to check that $\Sigma_N$ is a Bernoulli shift on  $p$ symbols, and that $\Phi_N: G \to X_N$ is a measurable (and topological) isomorphism.  Details of this construction are worked out in \cite{polMarkov}.
\end{example}

\section{Factoring through projections}\label{sec:factProj}
Let $G,H$ be compact topological groups.  For a transformation $T: G \to G$ we say that $T$ \emph{factors through} a surjective continuous group homomorphism $\phi: G \to H$ if there exists a transformation $T': H \to H$ such that the following diagram commutes
\[ \xymatrix{
G \ar[r]^T \ar[d]^\phi & G \ar[d]^\phi \\
H \ar[r]^{T'} & H } \]

Let us relate this to the situation of Lemma~\ref{lem:ilimProf}.
\begin{lemma}\label{lem:factorSigmaN}
Let $G$ be a profinite group, $\mu$ normalized Haar measure on $G$, and $T: G \to G$ a transformation on $G$.  Let $N \lhd_O G$ be such that $T$ factors through the quotient map $\pi_N: G \to G/N$.  Let $T'_N: G/N \to G/N$ denote the factor transformation.  Let $\Sigma_N = (X_N,\mu_N,\B_N,T_N)$ be as defined in Lemma~\ref{lem:ilimProf}.  Define
\[ \Sigma'_N = (G/N, \mu_{G/N}, \B_{G/N}, T'_N) \] where $\mu_{G/N}$ is Haar measure on the finite group $G/N$ (i.e. normalized counting measure), and $\B_{G/N}$ its $\sigma$-algebra (i.e. the power set of $G/N$).  Then, projection to the first coordinate $X_N \to G/N$ gives an isomorphism
\[ \Sigma_N \isom \Sigma'_N. \]
\end{lemma}
\begin{proof}
For $k \geq 0$, let $\varpi_k: X_N \to G/N$ be the projection to the $k^\text{th}$ coordinate.  Then, by the definition of $T_N$ and $T'_N$ we have the commutative diagram
\[ \xymatrix{  
     & G \ar '[d] [ddd]^{T^k} \ar[dr]^{\pi_N} \ar[dl]_{\Phi_N} &        \\
 X_A \ar[ddd]_{T_N^k} \ar[rr]_(.65){\varpi_0} &    &  G/N \ar[ddd]^{{T'_N}^k} \\
 & &               \\
     & G \ar[dr]^{\pi_N} \ar[dl]_{\Phi_N} &        \\
 X_A \ar[rr]^{\varpi_0}&    &  G/N   } \]
%\[ \xymatrix{ G \ar[r]^{T^k} \ar[d]^{\Phi_N} \ar@/_1pc/[dd]_{\pi_N} & G \ar[d]_{\Phi_N} \ar@/^1pc/[dd]^{\pi_N} \\ X_N \ar[d]^{\varpi_0} \ar[r]^{T_N^k} & X_N \ar[d]_{\varpi_0} \\ G/N \ar[r]^{{T'_N}^k} & G/N } \]
for each $k \geq 0$, where $T^k$, $T_N^k$, and ${T'_N}^k$ denote the $k$-fold composites of $T, T_N, T'_N$ respectively.

Note that for $x \in G/N$, 
\[ \mu_{N}\left(\varpi_0^{-1}(x)\right) = \mu\left(\Phi_N^{-1} \varpi_0^{-1}(x) \right) = \mu\left( \pi_N^{-1}(x) \right) = \mu_{G/N}(x) .\]  So, $\varpi_0$ is measure-preserving and $\Sigma'_N$ is a measurable factor of $\Sigma_N$.  Moreover, note that $\varpi_k = \varpi_0 \circ T_N^k = \varpi_0 \circ {T'_N}^k \varpi_0$; so each element of $X_N$ is uniquely determined by its first entry.  It follows that $\varpi_0^{-1}(\B_{G/N})=\B_N$.  Then,
\[ \Sigma_N' \isom (X_N,\mu_N,\varpi_0^{-1}(\B_{G/N}),T_N) = (X_N,\mu_N,\B_N,T_N) = \Sigma_N. \qedhere \]
\end{proof}

For $T: G \to G$, define the \emph{finite factor set of $T$} as
\[ \F(T) = \{ N \lhd_O G : T\text{ factors through }\pi_N:G \to G/N \}. \]  Note that each $\pi_N$ is a continuous surjective group homomorphism, thus measure-preserving with respect to Haar measure.

\begin{remark}
The notion of $\F(T)$ has another natural description.  Denote
\[ \F'(T) = \{ \pi \in \Hom_{\TopGp}(G,H)\text{ surjective}: H\text{ is a finite group},  T \text{ factors through }\pi \}{/\{\sim\}}, \]
where $\pi_1 \sim \pi_2$ if there exists an isomorphism $\im \pi_1 \isom \im \pi_2$ conjugating the two maps.

That is, $\F'(T)$ is the set of all finite group factors of $T: G \to G$.  The relationship  between $\F(T)$ and $\F'(T)$ is clear: for each $N \in \F(T)$ we have $G \to G/N \in \F'(T)$, and conversely for each $\pi \in \F'(T)$ we have $\ker \pi \in \F(T)$.
\end{remark}

\begin{defn}
For a profinite group $G$, we say that $T: G \to G$ is a \emph{{\psm}} if the cosets of $\F(T)$ form a base for the topology of $G$. 
\end{defn}

If $G$ is known to be second-countable, then Lemma~\ref{lem:ilimProf} and Lemma~\ref{lem:factorSigmaN} give us that
\[ \Sigma \eqdef (G,\mu,\B,T) \isom \ilim_{N \in \F(T)} \Sigma'_N, \] where $\Sigma'_N$ is in the sense of Lemma~\ref{lem:factorSigmaN}, and $\Sigma'_N$ is in particular a measurable dynamical system on a finite set.

We invite the reader to prove the following alternate characterization of the {\psm}s:
\begin{lemma}\label{lem:profSelfIlim}
Let $G$ be a profinite group.  Then a map $T: G \to G$ is a {\psm} if and only if there exists a directed set $(I,\leq)$ and an inverse system $\mathcal{D}: I \to \TopGp$ of finite groups such that
\[ G \isom \ilim_{i \in I} \mathcal{D}(i) \]
and $T$ factors through the projection $G \to \mathcal{D}(i) $ for each $i \in I$.  The inverse system may be assumed surjective.  In addition, instead of $T$ factoring through each projection, it suffices that for each $i \in I$ there exists a $j \in I$ with $i \leq j$ such that $T$ factors through the projection $G \to \D_j$.
\end{lemma}

\begin{example}
Let $G = \ZZ_p$.  We note that
\[ \ZZ_p \isom \ilim_{k \in \NN} \ZZ/p^k \ZZ. \]

The open normal subgroups of $\ZZ_p$ are all of the form $p^k \ZZ_p$.  Then, we see that $T: \ZZ_p \to \ZZ_p$ is a {\psm} if and only if 
there is an infinite subset $I$ of $\NN$ such that $k \in I$ and $|x-y|_p \leq p^{-k}$ implies that $|Tx-Ty|_p \leq p^{-k}$.  In particular, this holds for all maps satisfying $|Tx-Ty|_p \leq |x-y|_p$ (i.e., the $1$-Lipschitz maps). In this context,  our notion of {\psm}s may be viewed as generalizing the notions of (asymptotically) compatible maps found in \cite{Anashin}, and our Proposition~\ref{prop:profMP} generalizes Lemma~4.5 of \cite{BrykSilva}.
\end{example}

\begin{example}
Let $G = \ZZ_p \times \ZZ_p$.  We note that
\[ G \isom \ilim_{k_1,k_2\in \NN} \ZZ/p^{k_1}\ZZ \times \ZZ/p^{k_2}\ZZ. \]

Let $T$ be given by multiplication by an element of $\GL_2(\ZZ_p)$.  Given $k_1,k_2 \in \NN$ it need not be the case that $T$ factors through the projection to $\ZZ/p^{k_1}\ZZ \times \ZZ/p^{k_2} \ZZ$.  However, $T$ does factor through the projection for $k_1=k_2$.  The kernels of these projections form a base for the neighborhoods of $\idgp \in G$, so $T$ is a {\psm}.
\end{example}

\begin{lemma}
Let $G$ be a profinite group and $T: G \to G$ a {\psm}.  Then, $T$ is continuous.
\end{lemma}
\begin{proof}
Say $T$ factors through $\pi_N: G \to G/N$ as $T_N$ for each $N \in \F(T)$.

For $N \in \F(T)$ and $h \in G/N$, then
\[ T^{-1}(\pi_N^{-1}(h)) = \pi_N^{-1}(T_N^{-1}(h)) = \bigcup_{h' \in T_N^{-1}(h)} \pi_N^{-1}(h'). \]

As the sets
\[ \{ \pi_N^{-1}(h) : N \in \F(T), h \in G/N \} \]
are precisely the cosets of the elements of $\F(T)$ they form a base for the topology on $G$.
As $T^{-1}$ takes each set in this base to an open set, continuity of $T$ follows.
\end{proof}

\begin{lemma}\label{lem:nonSingSurj}
Let $G$ be a compact Hausdorff topological group, $\mu$ normalized Haar measure on $G$, and $T: G \to G$ continuous.  If $T$ is nonsingular with respect to $\mu$, then $T$ is surjective.
\end{lemma}
\begin{proof}
As $T$ is continuous, $T(G)$ is the continuous image of a compact set, thus compact and so closed in the Hausdorff space $G$.

But,
\[ \mu\left(T^{-1}\left(G \setminus T(G)\right) \right) = \mu\left(\emptyset\right) = 0, \]
and by nonsingularity
\[ \mu\left(G \setminus T(G)\right) = 0. \]

Note that $\mu$ is positive on non-empty open sets, so this implies that  $G \setminus T(G)$ does not contain a non-empty open set and hence that $T(G)$ is dense in $G$.  As $T(G)$ is closed in $G$, this implies $T(G)=G$.  So, $T$ is surjective.
\end{proof}

\begin{prop}\label{prop:profMP}
Let $G$ be a second-countable profinite group, $\mu$ normalized Haar measure on $G$, and $T: G \to G$ a {\psm}.  Let $\F \subseteq \F(T)$ be a base for the neighborhoods of $\idgp \in G$.  For each $N \in \F(T)$ let $T_N$ denote the induced map $G/N \to G/N$.  Then, the following are equivalent:
\begin{enumerate}
\item $T_N$ is bijective on $G/N$ for all $N \in \F$;
\item $T_N$ is nonsingular with respect to $\mu_{G/N}$ for all $N \in \F$;
\item $T_N$ is measure-preserving with respect to $\mu_{G/N}$ for all $N \in \F$;
\item $T$ is measure-preserving with respect to $\mu$;
\item $T$ is nonsingular with respect to $\mu$;
\item $T$ is surjective.
\end{enumerate}
\end{prop}
\begin{proof}
We prove the following implications:
\[ \xymatrix@C=.5cm{
 \text{(i)} \ar@{<=>}[r] & \text{(ii)} \ar@{<=>}[r]& \text{(iii)} \ar@{<=>}[r] & \text{(iv)} \ar@{=>}[r] & \text{(v)} \ar@{=>}[r] & \text{(vi)} \ar@{=>}[r] & \text{(i)} } \]
The implications (i)$\Leftrightarrow$(ii)$\Leftrightarrow$(iii) follow as $G/N$ is finite and $\mu_{G/N}$ is counting measure.  The implications (iii)$\Leftrightarrow$(iv) follow by Lemma~\ref{lem:ilimProf} and Lemma~\ref{lem:factorSigmaN}.  The implication (iv)$\Rightarrow$(v) is true by definition.  The implication (v)$\Rightarrow$(vi) follows by Lemma~\ref{lem:nonSingSurj}.  Finally, (vi) implies that each $T_N$ is surjective; as a surjective map of a finite set to itself is bijective, we have (vi)$\Rightarrow$(i).
\end{proof}

\begin{lemma}\label{lem:profMPMetric}
Let $G$ be a second-countable profinite group, $\mu$ normalized Haar measure on $G$, and $T: G \to G$ a transformation.  Then, $T$ is a measure-preserving {\psm} if and only if there exists a metric $d: G^2 \to \RR_{\geq 0}$ on $G$ such that the following conditions hold:
\begin{enumerate}
\item $d$ induces the usual topology on $G$;
\item $d$ is left translation invariant in the sense that $d(gx,gy)=d(x,y)$ for all $x,y,g \in G$;
\item $T$ is an isometry with respect to $d$.
\item The set of open-subgroups of $G$ which are (closed) balls with respect to $d$, i.e. \[ \left\{ N \lhd_O G : N = \{ x \in G : d(\idgp, x) \leq r_N \}\text{ for some $r_N > 0$} \right\} \] is a base for the neighborhoods of $\idgp \in G$. %There exists a collection $\F \subseteq \{ N \lhd_O G \}$, such that $\F$ is a base for the neighborhoods of $\idgp \in G$, and such that for each $N \in \F$ there exists an $r_N \in \RR_{\geq 0}$ such that for each $g \in G$ it is the case that $N = \{ x \in G : d(\idgp, x) \leq r_N \}$ (that is, the elements of $\F$ are balls with respect to $d$);
\end{enumerate}
\end{lemma}
\begin{proof}
{\noindent}{\bf $\Rightarrow$: }\\
By the proof of Lemma~\ref{lem:ilimProf} we note that $\F(T)$ is countable and that the translates of the elements of $\F(T)$ give a countable base for the topology on $G$.  Say $\F(T) = \{ N'_1, N'_2,N'_3,\ldots \}$.  Set $N_1 = N'_1$, and for $k > 1$ let $N_k \in \F(T)$ be such that $N_k \subseteq N_{k-1} \cap N'_k$. Note that $N_{k-1} \cap N'_k$ is open and contains $\idgp$ for each $k > 1$, so such an $N_k$ must exist.  Then, set $\F = \{ N_1, N_2, \ldots \}$.  Note that $\F$ is countable, nested, and forms a base for the neighborhoods of $\idgp \in G$.  

For $N \lhd_O G$, let $\pi_N: G \to G/N$ be the quotient map.  Then, we may define $d: G^2 \to \RR_{\geq 0}$ for $x,y \in G$ by
\[ d(x,y) = 2^{-\ell} \text{ where $\ell = \min \{ k : \pi_{N_k}(x) = \pi_{N_k}(y) \}$}, \]
and $d(x,y)=0$ if $\pi_{N_k}(x)=\pi_{N_k}(y)$ for all $k \geq 0$.  

We claim that $d$ is a metric, and that it moreover satisfies the conditions in the lemma:
\begin{itemize}
\item  It is clear by construction that $d$ is symmetric and non-negative.  Note that
\[ \bigcap_{k \geq 0} N_k = \{ \idgp \}, \] so $d(x,y) = 0 \Leftrightarrow x = y$.  Moreover, $d(x,y) \leq 2^{-k}$ and $d(y,z) \leq 2^{-k}$ implies $d(x,z) \leq 2^{-k}$, so $d$ satisfies the strong triangle inequality.  So, we see that $d$ is indeed a metric.
\item The set of balls with respect to $d$ is precisely $\F$ and the emptyset.  So, $d$ satisfies condition (iv) of the Lemma, and moreover it induces the same topology as $\F$ and so satisfies (i).
\item As $\pi_{N_k}$ is a homomorphism for each $k \geq 0$, we see immediately that $\pi_{N_k}(x) = \pi_{N_k}(y) \Leftrightarrow \pi_{N_k}(gx) = \pi_{N_k}(gy) \Leftrightarrow \pi_{N_k}(xg)=\pi_{N_k}(yg)$ for all $x,y,g \in G$.  So, $d$ is (left and right) translation invariant, and satisfies (ii).
\item For $N \in \F \subseteq \F(T)$ we have that $T_N: G/N \to G/N$ is a bijection by Proposition~\ref{prop:profMP}.  So, $\pi_{N_k}(x) = \pi_{N_k}(y) \Leftrightarrow \pi_{N_k}(T(x)) = \pi_{N_k}(T(y))$ for all $k \geq 0$.  So, we see that $T$ is an isometry with respect to $d$, hence condition (iii).
\end{itemize}

\medskip
{\noindent}{\bf $\Leftarrow$: }\\
Let $\F$ be the collection in (iv).  For each $N \in \F$, let $r_N > 0$ be as in the definition of $\F$.

Using the fact that $N = \{ x \in G : d(\idgp, x)=d(x,\idgp) \leq r_N \}$ and the fact that $d$ is left translation invariant we confirm that 
\[ \pi_N(x) = \pi_N(y) \Leftrightarrow x^{-1} y \in N \Leftrightarrow d(x,y) = d(x^{-1} y,\idgp) \leq r_N. \]    Then, the fact that $T$ is an isometry with respect to $d$ implies that $\pi_N(x) = \pi_N(y) \Leftrightarrow \pi_N(T(x)) = \pi_N(T(y))$.  So, $T$ induces a well-defined injective, hence bijective as $G/N$ is finite, map $T_N: G/N \to G/N$.

As this holds for arbitrary $N \in \F$, we have that $\F \subseteq \F(T)$ is a base for the neighborhoods of $\idgp \in G$ with $T_N$ bijective on $G/N$ for all $N \in \F$.  This implies immediately that $T$ is a {\psm}, and by Proposition~\ref{prop:profMP} that $T$ is measure-preserving.
\end{proof}

\begin{prop}\label{prop:profErg}
Let $G$ be a second-countable profinite group, $\mu$ normalized Haar measure on $G$, and $T: G \to G$ a {\psm}.  Let $\F \subseteq \mathcal{F}(T)$ be a base for the neighborhoods of $\idgp \in G$.  For each $N \in \mathcal{F}(T)$ let $T_N$ denote the induced map $G/N \to G/N$.  Then, the following are equivalent:
\begin{enumerate}
\item $T$ is measure-preserving and ergodic with respect to $\mu$;
\item $T_N$ is measure-preserving and ergodic with respect to $\mu_{G/N}$ for all $N \in \F$;
\item $T_N$ is minimal for all $N \in \F$.
\end{enumerate}
By Proposition~\ref{prop:profMP}, we may replace ``measure-preserving'' with ``nonsingular'' in one or both of the above occurrences.
\end{prop}
\begin{proof}
The equivalence (i)$\Leftrightarrow$(ii) follows by Lemma~\ref{lem:ilimProf} and Lemma~\ref{lem:factorSigmaN}.  The equivalence (ii)$\Leftrightarrow$(iii) holds as each $G/N$ is finite with $\mu_{G/N}$ the normalized counting measure.
\end{proof}

\begin{prop}\label{prop:profWeakMix}
Let $G$ be a compact Hausdorff topological group, $\mu$ normalized Haar measure on $G$, and $T: G \to G$ a transformation.  Say $\F(T) \supsetneq \{G\}$.  Then, $T$ is not totally ergodic .  In particular, if $T$ is a {\psm} then it is not totally ergodic unless $|G|=1$.
\end{prop}
\begin{proof}
Let $N \in \F(T) \setminus \{ G \}$.  Then, $1 < |G/N|  < \infty$, and $T: G \to G$ factors through $G/N$ as
\[ \xymatrix{
G \ar[r]^T \ar[d]^{\pi_N} & G \ar[d]^{\pi_N} \\
G/N \ar[r]^{T_N} & G/N } \]

If $T$ is ergodic then $T_N$ is ergodic, hence minimal.  In particular for $h \in G/N$ we have that $T^\ell(h)=h$ if and only if $|G/N| \divides \ell$.  Then, $T^{|G/N|}$ factors through the projection as $T_{N}^{|G/N|}$; but this is just the identity map on $G/N$.  So, $T$ is not totally ergodic.

If $T$ is a {\psm}, then \[ G \isom \ilim_{N \in \F(T)} G/N \] by Proposition~\ref{prop:propProfGp}.  In particular $\F(T) = \{ G \}$ implies $|G| = 1$.
\end{proof}

%Note that an isometry (or an equicontinuous map in general) can not be totally ergodic.  So, the final sentence of Proposition~\ref{prop:profWeakMix} also follows immediately from Lemma~\ref{lem:profMPMetric}.
\begin{remark}
Recall also that weakly mixing implies totally ergodic.  So, the above also gives negative results for weak mixing.
\end{remark}

Now,  the results of the propositions yield the proof of Theorem~\ref{thm:profDynamics}.  

\section{Homomorphisms}\label{sec:hom}
We begin by recalling a result on when a continuous group endomorphism is measure-preserving:
\begin{lemma}\label{lem:homMP}
Let $G$ be a compact Hausdorff topological group, $\mu$ normalized Haar measure on $G$, and $T: G \to G$ a homomorphism of topological groups.  Then, the following are equivalent
\begin{enumerate}
\item $T$ is nonsingular with respect to $\mu$;
\item $T$ is surjective;
\item $T$ is measure-preserving with respect to $\mu$.
\end{enumerate}
\end{lemma}
\begin{proof}
The assertion (i)$\Rightarrow$(ii) follows from Proposition~\ref{prop:profMP}.  The assertion (ii)$\Rightarrow$(iii) is true as $\mu \circ T^{-1}$ can be shown to be regular, translation invariant, and normalized.  The assertion (iii)$\Rightarrow$(i) is true by definition.
\end{proof}

Now, in the case of continuous group endomorphisms, we may give an alternate characterization of the collection $\F(T)$ in the definition of a {\psm}:
\begin{lemma}\label{lem:homFactors}
Let $G$ be a compact Hausdorff topological group and $T: G \to G$ a homomorphism of topological groups.  Then, \[ \F(T) = \{ N \lhd_O G : N \subseteq T^{-1}(N) \}. \]  If $T$ is surjective then in fact \[ \F(T) = \{ N \lhd_O G : N = T^{-1}(N) \}. \]
\end{lemma}
\begin{proof}
Note that for $N \lhd_O G$, any $T_N$ making the following diagram commute must be a group homomorphism
\[ \xymatrix{
G \ar[r]^T \ar[d]^{\pi_N} & G \ar[d]^{\pi_N} \\
G/N \ar[r]^{T_N} & G/N } \]

Furthermore, such a $T'$ exists if and only if $N = \ker \pi_N \subseteq \ker \pi_N \circ T = T^{-1}(N)$.  If $T$ is in addition surjective, then by Lemma~\ref{lem:homMP} it is measure-preserving.  Then, $\mu(N) = \mu(T^{-1}(N))$ and so $\mu(T^{-1}(N) \setminus N) = 0$; as $T^{-1}(N) \setminus N$ is open, this implies that it is empty and so $T^{-1}(N)=N$.

So, \[ \F(T) = \{ N \lhd_O G : T(N) \subseteq N \},\] and if $T$ is in addition surjective then we may replace the constraint by $T(N) = N$.
\end{proof}

\begin{remark}\label{rem:profGpNotErg}
Note that if $\F(T) \neq \{ G \}$ then $T$ is not ergodic.  This follows because the factor transformation would be a group homomorphism on a finite group, which can not be ergodic (for it maps $\idgp$ to itself).
\end{remark}

For many profinite groups, the following criterion suffices to show that all group endormorphisms are {\psm}s:
\begin{prop}\label{prop:homNormal1}
Let $G$ be a profinite group such that $G$ has finitely many open normal subgroups of each finite index.  If $T: G \to G$ is a (Haar) nonsingular homomorphism of topological groups (i.e. a surjective continuous group homomorphism), then $T$ is a {\psm}.

In particular, if $G$ has a finitely-generated dense subgroup then the any such $T$ is a {\psm}.
\end{prop}
\begin{proof}
Say $N \lhd_O G$.  Then, $T^{-1}(N) \lhd_O G$.  Taking measures and noting that $T$ is measure-preserving with respect to Haar measure by Lemma~\ref{lem:homMP} we observe that
\[ 1/[G:N]=\mu(N) = \mu(T^{-1}(N)) = 1/[G:\mu(T^{-1}(N)). \]

Now, for $N \lhd_O G$ consider the collection
\[ \{ T^{-k}(N) : k \geq 0 \}. \]  Each element of the this collection must be an open normal subgroup of the same index in $G$, so the collection must be finite by hypothesis.  Set \[ N' = \bigcap_{k \geq 0} T^{-k}(N), \] where the intersection is over finitely many distinct sets; so $N' \lhd_O G$.  Note that $N \cap T^{-1}(N') = N'$, so $N' \subseteq T^{-1}(N')$ and $N' \in \F(T)$ by Lemma~\ref{lem:homFactors}.  Moreover, $N' \subseteq N$ and $N$ may be written as a union of cosets of $N'$.  As this holds for arbitrary $N \lhd_O G$, we see that $\F(T)$ forms a base for the neighborhoods of $\idgp \in G$, and $T$ is a {\psm}.

By \cite[Lemma~4.1.2]{Wilson}, if $G$ has a finitely-generated dense subgroup then $G$ has finitely many open normal subgroups of a given index, and the final assertion of the proposition follows.
\end{proof}

%
%\begin{remark}
%The lemma is of course not the only method of proving that a topological group has finitely many open normal subgroups of a given index.  One classical example comes from Galois theory of the $p$-adics:
%Let $\ol{\QQ_p}$ be an algebraic closure of $\QQ_p$, and let $G = \Gal(\ol{\QQ_p}/\QQ_p)$ be the Galois group of this extension.  By the Galois correspondence of infinite Galois theory, we have that closed subgroups of index $n < \infty$ in $G$ correspond bijectively to Galois extensions of degree $n$ over $\QQ_p$.  It is a standard application of Krasner's Lemma and the classification of unramified and totally ramified extensions to show that there are finitely many extensions of $\QQ_p$ of a given degree; the interested reader may consult \cite{Koblitz}, particularly Ch. III \S 3 and Ch. III Exercise~4, and note that the result follows by noting that the space of Eisenstein polynomials of a given degree is compact (for $\pp$ is compact) and using Exercise~4.  It follows that $G$ has finitely many closed subgroups of any given index.
%\end{remark}

We may apply the Propositon to several groups of interest:
\begin{corollary}
Let $G = \prod_{i=1}^{g} \ZZ_{p_i}^{e_i}$ with the $p_i$ rational primes and $e_i \in \NN$.  Then, any continuous homomorphism $T: G \to G$ is a {\psm} and is not ergodic.
\end{corollary}
\begin{proof}
We note that $G$ contains a dense finitely-generated subgroup
\[ \prod_{i=1}^{g} \ZZ^{e_i}. \]

Then, $G$ has finitely many open normal subgroups of a given index, and in particular for each open normal subgroup $N \lhd_O G$ we have that $\{ T^{-k}(N) \}$ must be finite (for each element of this set has index equal to the index of $N$).  Applying Proposition~\ref{prop:homNormal1} proves that $T$ is a {\psm}, and applying Remark~\ref{rem:profGpNotErg} yields that $T$ is not ergodic.
\end{proof}

\begin{corollary}
Let $G = \ZZ_p^{k}$.  Then, the nonsingular continuous homomorphisms $T: G \to G$ are given by multiplication by elements of $\GL_k(\ZZ_p)$.  Any such homomorphism is a {\psm} and is not ergodic.
\end{corollary}
\begin{proof}
We note that $\ZZ^k$ is dense in $G$, and so a continuous homomorphism is defined by its values on a basis for $\ZZ^k$.  In particular, this implies that any continuous homomorphism must be given by multiplication by some $T \in \Mat_{k \times k}(\ZZ_p)$.  By Proposition~\ref{prop:profMP} we must have $T$ surjective.  In particular, the image of $T$ must contain the generators for $\ZZ_p^k$, so there must exist a $S \in \Mat_{k \times k}(\ZZ_p)$ such that $TS=\id_{k \times k} \in \Mat_{k \times k}(\ZZ_p)$.  Then, $T \in \GL_k(\ZZ_p)$ (and of course, the converse holds by reversing this logic).  Now, the previous corollary gives that this map must be a {\psm} and is not ergodic.
\end{proof}

\begin{remark}
In this context we mention that Juzvinski{\u\i}  \cite{Juzvinskii} 
showed that ergodic group endomorphisms have completely positive 
entropy and    Lind proves in \cite{Lind} that ergodic automorphisms of 
compact metrizable groups are measurably isomorphic to Bernoulli 
shifts.
\end{remark}

\section{Products}\label{sec:prod}
\begin{lemma}\label{lem:prodProf}
Let $A$ be an index set.  For each $\alpha \in A$ let $G_\alpha$ be a profinite group and $T_\alpha: G_\alpha \to G_\alpha$ a {\psm}.  Then, \[ G = \prod_{\alpha \in A} G_\alpha \] is a profinite group, and \[ T = \prod_{\alpha \in A} T_\alpha \] is a {\psm} on $G$.
\end{lemma}
\begin{proof}
Note that for each $\alpha \in A$ we have that $\F(T_\alpha)$ is a base for the neighborhoods of $\idgp \in G_\alpha$.  Then, the collection
\[\F = \{ \prod_{\alpha \in A} N_\alpha : N_\alpha \in \F(T_\alpha), N_\alpha = G_\alpha\text{ for all but finitely many } \alpha \in A \} \]
forms a base for the neighborhoods of $\idgp \in  G$.  Moreover, observe that each element of $\F$ is a normal subgroup of $G$.

We claim that
\[ G \isom \ilim_{N \in \F} G/N. \]
Indeed, the natural projections induce a homomorphism
\[ \xymatrix{ \phi: G \ar[r]^{\phi} & \ilim_{N \in F} G/N }. \]
Observe that $\phi$ is injective as $G$ is Hausdorff.  Moreover, $\phi$ is continuous and $G$ compact (by Tychonoff's Theorem), so the image of $\phi$ is closed; but the image of $\phi$ is also dense in the codomain.  So, $\phi$ is surjective.  Then, $\phi$ is a continuous bijection with compact domain, so a homeomorphism, and $G$ is indeed profinite.

Now, note that for any $N \in \F$, $T$ factors through the projection $G \to G/N$ as the product of the factor transformations in each coordinate.  So, $T$ is a {\psm}.
\end{proof}

\begin{lemma}\label{lem:crtPlus}
Let $S_k$ be a finite non-empty set and $T_k: S_k \to S_k$ a transformation for $k=1,\ldots,n$.  Let
\[ S = \prod_{k=1}^{n} S_k, \qquad T = \prod_{k=1}^{n} T_k. \]  Then, $T$ is minimal on $S$ if and only if each $T_k$ is minimal on $S_k$ and the $|S_k|$ are pairwise coprime.
\end{lemma}
\begin{proof}
Note that the general case follows from $n=2$ case by induction.  So, we may assume $n=2$.

We have that $T$ minimal implies $T_1,T_2$ minimal.  By the minimality of $T_k$, each point of $S_k$ must have full orbit.  So we have $T_k^\ell(x)=x$ if and only if $|S_k| \divides \ell$. Let $\ell = |S_1| |S_2| / (|S_1|,|S_2|)$ be the least common multiple of $|S_1|, |S_2|$.  Then, 
\[ T^\ell(s_1,s_2) = (T_1^\ell(s_1),T_2^\ell(s_2))=(s_1,s_2). \]  So, $T$ minimal requires $(|S_1|,|S_2|)=1$, that is that the cardinalities be coprime.

Conversely, say $(|S_1|,|S_2|)=1$.  In particular, given $s_k \in T_k$, $\ell_k \in \NN$ for $k=1,2$, the Chinese Remainder Theorem gives us a $\ell \in \NN$ such that $\ell \equiv \ell_k \pmod{|S_k|}$ for $k=1,2$.  Then,
\[ T^\ell(s_1,s_2) = (T_1^\ell(s_1),T_2^\ell(s_2))=(s_1^{\ell_1},s_2^{\ell_2}). \]  Then, $T_1, T_2$ minimal implies $T$ minimal.
\end{proof}

Then:
\begin{theorem}\label{thm:prodProfErg}
Let $A, G_\alpha, T_\alpha, G, T$ be as in Lemma~\ref{lem:prodProf}.  Moreover, assume each $G_\alpha$ is second-countable and $A$ is countable. Then, $G$ is second-countable and
\begin{enumerate}
\item $T$ is nonsingular if and only if $T_\alpha$ is nonsingular for each $\alpha \in A$.
\item Denote \[ D_\alpha = \{ |G_\alpha/N_\alpha| : N_\alpha \in \F(T_\alpha) \}. \]  Then, $T$ is ergodic if and only if $T_\alpha$ is ergodic for each $\alpha \in A$ and for all $\alpha,\beta \in A$ distinct and all $n \in D_\alpha, m \in D_\beta$ we have $(n,m)=1$.
\end{enumerate}
\end{theorem}
\begin{proof}
For each $\alpha \in A$ let $C_\alpha$ be a countable base for $G_\alpha$.  We may assume without loss of generality that $G_\alpha \in C_\alpha$ for each $\alpha \in A$.  Then the set
\[ \{ \prod_{\alpha \in A} S_\alpha : S_\alpha\in C_\alpha, S\alpha = G_\alpha \text{ for all but finitely many $\alpha \in A$} \} \]
is a countable base for $G$.  So, $G$ is second-countable.

Note that $T$ and each $T_\alpha$ are {\psm}s.  So, by Proposition~\ref{prop:profMP}, they are nonsingular if and only if they are surjective.  Now, the product of a set of maps is surjective if and only if each map is surjective.  The first claim follows.

Applying Proposition~\ref{prop:profErg} to each $T_\alpha$ we see that $T_\alpha$ is ergodic if and only if each of the factor transformations $\{ T_\alpha^{N_\alpha} : N_\alpha \in \F(T_\alpha) \}$ is minimal.  

For $N_\alpha \in \F(T_\alpha)$, let $T_\alpha^{N_\alpha}$ denote the map making the following diagram commute \[ \xymatrix{ G_\alpha \ar[r]^{T_\alpha} \ar[d] & G_\alpha \ar[d] \\ G_\alpha/N_\alpha \ar[r]^{T_\alpha^{N_\alpha}} & G_\alpha/N_\alpha} \]

We note that $\F(T_\alpha)$ is a base for the open sets containing  $\idgp \in G_\alpha$, and so,
\[ \{ \prod_{\alpha \in A} N_\alpha : N_\alpha \in \F(T_\alpha), N_\alpha = G_\alpha\text{ for all but finitely many } \alpha \in A \} \]
is a base for the open sets containing $\idgp \in G$.  Given \[ N = \prod_{\alpha \in A} N_\alpha \] in this base, we have that $T$ factors through the projection $G \to G/N$ as 
\[ T_N = \prod_{\alpha \in A} T_{\alpha}^{N_\alpha}. \]  

Applying Proposition~\ref{prop:profErg}, we see that $T$ is ergodic if and only if each of these factor transformations is minimal on the finite quotient
\[ G/N = \prod_{\alpha \in A} G_{\alpha}/N_{\alpha} \isom \prod_{N_\alpha \neq G_\alpha} G_{\alpha}/N_\alpha. \]  Dropping trivial factors and applying Lemma~\ref{lem:crtPlus}, we have that $T$ is ergodic if and only if each $T_{\alpha}^{N_\alpha}$ is ergodic for all $\alpha \in A$ and $N_\alpha \in \F(T_\alpha)$ and the elements of the $D_\alpha$ are pairwise co-prime (for different subscripts).  Applying Proposition~\ref{prop:profErg} to the $T_\alpha$, this yields our desired result.
\end{proof}

\begin{remark}
As a consequence, we get an alternate proof that no {\psm}s are weakly mixing.
\end{remark}

\begin{corollary}
Let $T_p: \ZZ_p \to \ZZ_p$ be an ergodic {\psm} for each rational prime $p$.  Then, the map $T = \prod_{p} T_p$ on $G = \prod_p \ZZ_p$ is an ergodic {\psm}.
\end{corollary}
\begin{proof}
Follows immediately by Theorem~\ref{thm:prodProfErg} after noting that $\ZZ_p$ has quotients of $p$-power orders.
\end{proof}

\begin{corollary}
The maps $x \mapsto x \pm 1$ on 
\[ \oh{\ZZ} = \ilim_{n,\divides} \ZZ/n\ZZ \isom \prod_{p} \ZZ_p \]
is ergodic.
\end{corollary}
\begin{proof}
Each maps factors through all the projections and so is a {\psm}. Note that the maps $x \mapsto x \pm 1$ are certainly minimal on $\ZZ/n\ZZ$ for each $n > 0$.  In light of Proposition~\ref{prop:profErg} this gives a direct proof that the induced map on $\oh{\ZZ}$  is ergodic.  Alternatively, we may use Proposition~\ref{prop:profErg} to show that the induced map on $\ZZ_p$ is ergodic for each $p$, and then use the previous corollary.

Also, observe that for $n>2$, the maps $x \mapsto -x \pm 1$ are \emph{not} minimal on $\ZZ/n\ZZ$ [$1-0=1$, $1-1=0$; $-1-0=-1$, $-1-(-1)=0$].
\end{proof}

\begin{remark}
Let $K = \FF_p$ be the finite field of $p$ elements.  Let $L$  be an algebraic closure of $K$.  Then, 
\[ G = \Gal(L/K) \isom \oh{\ZZ}. \]   The $p^\text{th}$ power map (the ``Frobenius automorphism''), denoted $\Frob \in G$, generates a dense cyclic subgroup of $G$.  Indeed, the map $\ZZ \to G$ given by $n \mapsto \Frob^n$ induces the above isomorphism.  So, the map $x \mapsto x+1$ on $\oh{\ZZ}$ may be reinterpreted as the map on $G$ given by $\sigma \mapsto \sigma \circ \Frob$.

Alternatively, we could let $K = \CC(t)$, and $L=\CC(t,t^{1/2},t^{1/3},t^{1/4},\ldots,t^{1/n},\ldots)$.  Then, \[ G = \Gal(L/K) \isom \oh{\ZZ}. \]  Take $\tau \in G$ defined by
\[ \tau t^{1/m} = e^{2\pi \ii/m} t^{1/m}. \]  Then, $\tau$ generates a dense cyclic subgroup of $G$, and the map $\sigma \mapsto \sigma \circ \tau$ is the equivalent of $x \mapsto x+1$.
\end{remark}

%Bibliography
%\bibliographystyle{amsalpha}
%\bibliography{ProfiniteKLPS}

\end{document}